\documentclass[a4paper,11pt]{article}
\usepackage[T1]{fontenc}
\usepackage[cp1250]{inputenc}
\usepackage{lmodern}

\usepackage[english]{babel}
\usepackage{latexsym}
\usepackage{amsmath}
\usepackage{indentfirst}
\usepackage{amsthm}
\usepackage{amssymb}
\usepackage{amsfonts}
\usepackage{stackrel}
\usepackage{dsfont}
\usepackage{tikz}
\usepackage{slashbox}
\usepackage[mathscr]{eucal}
\usepackage{authblk}

\newtheorem{thm}{Theorem}
\newtheorem{lemma}[thm]{Lemma}
\newtheorem{cor}[thm]{Corollary}

\newenvironment{pro}{\begin{flushleft} \textbf{Proof}\\* \end{flushleft}}{\hfill\(\blacksquare\) \\ }

\newcommand{\reals}{\mathbb{R}}

\newcommand{\complex}{\mathbb{C}}

\newcommand{\eps}{\varepsilon}

\newcommand{\Ffamily}{\mathcal{F}}

\newcommand{\Ufamily}{\mathcal{U}}
\newcommand{\Vfamily}{\mathcal{V}}
\newcommand{\Pfamily}{\mathcal{P}}

\newcommand{\Wall}{\text{Wall}}
\newcommand{\Zero}{\mathcal{Z}}

\newcommand{\Addresses}{{
  \bigskip
  \footnotesize

  Krukowski M. (corresponding author), \textsc{Technical Univeristy of \L\'od\'z, \ Institute of Mathematics, \ W\'ol\-cza\'n\-ska 215, \
90-924 \ \L\'od\'z, \ Poland}\par\nopagebreak
  \textit{E-mail address} : \texttt{krukowski.mateusz13@gmail.com}

  \medskip
}}

\begin{document}
\title{Arzel\`a-Ascoli theorem via Wallman compactification}
\author{Mateusz Krukowski}
\affil{Technical University of \L\'od\'z, Institute of Mathematics, \\ W\'ol\-cza\'n\-ska 215, \
90-924 \ \L\'od\'z, \ Poland}
\maketitle

\begin{abstract}
In the paper, we recall the Wallman compactification of a Tychonoff space $T$ (denoted by $\Wall(T)$) and the contribution made by Gillman and Jerison. Motivated by the Gelfand-Naimark theorem, we investigate the homeomorphism between $C^b(T)$ and $C(\Wall(T))$. Along the way, we attempt to justify the advantages of Wallman compactification over other manifestations of Stone-\v{C}ech compactification. The main result of the paper is a new form of Arzel\`a-Ascoli theorem, which introduces the concept of equicontinuity along $\omega$-ultrafilters.
\end{abstract}


\section{Introduction}

The classical Arzel\`a-Ascoli theorem (\cite{Munkres} on page 278) plays an important role in functional analysis. It characterizes the relatively compact subsets of the space of complex-valued continuous functions $C(X)$, with $X$ a compact space, as those which are equibounded and equicontinuous. The space $C(X)$ is given the standard norm 

$$\|f\| := \sup_{x\in X} |f(x)|$$

The theorem admits numerous generalizations. A version for locally compact space $X$ and metric space $Y$ can be found in \cite{Munkres} on page 290. The topology on $C(X,Y)$ is the topology of uniform convergence on compacta.

In \cite{Przeradzki}, Bogdan Przeradzki studied the existence of bounded solutions to the equations $x' = A(t)x + r(x,t)$, where $A$ is a continuous function taking values in the space of bounded linear operators in a Hilbert space and $r$ is a nonlinear continuous mapping. He came up with a characterization of relatively compact subsets of the space of bounded and continuous functions $C^b(\reals,E)$, where $E$ is a Banach space. In addition to pointwise relative compactness and equicontinuity, the following condition was introduced:

\begin{description}
	\item[\hspace{0.4cm} (P)] For any $\eps > 0$, there exist $T> 0$ and $\delta > 0$ such that if $\|x(T) - y(T)\| \leq \delta$ then $\|x(t) - y(t)\| \leq \eps$ for $t \geq T$ and if $\|x(-T) - y(-T)\|\leq \delta$ then $\|x(t) - y(t)\|\leq \eps$ for $t \leq -T$ where $x$ and $y$ are arbitrary functions in $\Ffamily$. 
\end{description}

This idea was further studied by Robert Sta\'nczy in \cite{Stanczy}, while investigating the existence of solutions to Hammerstein equations in the space of bounded and continuous functions. The results were applied to Wiener-Hopf equations and to ODE's. In \cite{KrukowskiPrzeradzki} the author, together with Bogdan Przeradzki recast \textbf{(P)} in a setting where $X$ is $\sigma$-locally compact Hausdorff space and $Y$ is a complete metric space. In this paper, we recall the obtained result (as theorem \ref{AAforXY}) without proof.

Bogdan Przeradzki suggested that Arzel\`a-Ascoli theorem could be described in terms of ultrafilters. The motivation comes from Gelfand-Naimark theorem in \cite{GelfandNaimark}, which states that $C^b(T)$ and $C(\beta T)$ are $^{\ast}$-isomorphic as $C^{\ast}$-algebras (so in particular homeomorphic), where $\beta T$ is the Stone-\v{C}ech compactification. The idea is to write the classical Arzel\`a-Ascoli theorem for $C(\beta T)$ and to interpret it from the perspective of $C^b(T)$. To this end, we will need an explicit form of the homeomorphism between these spaces.

We will use Wallman ultrafilter construction rather than the original approach by Stone or \v{C}ech. The reason is that Wallman topology will be more convenient to work with than the weak$^{\ast}$ topology of $\beta T$. We attempt to explain the advantages of this topology at the beginning of section \ref{sectionWallman}. Moreover, Wallman's original construction in \cite{Wallman} has been improved by Gillman and Jerison in \cite{GillmanJerison} and even further by Frink in \cite{Frink} or Steiner in \cite{Steiner}. Although, we will not need full generality, the crucial parts of the construction are briefly summarized in this section.

The main results are described in section \ref{sectionmainpart}. Theorem \ref{homeobc} describes the homeomorphism between $C^b(T)$ and $C(\Wall(T))$, as anticipated by the Gelfand-Naimark theorem. The culminating point is theorem \ref{AA}, which characterizes relatively compact subsets of $C^b(T)$ in terms of $\omega$-ultrafilters. The theorem introduces the concept of equicontinuity along $\omega$-ultrafilters, which accounts for classical equicontinuity and $C^b(X,Y)$-extension property in theorem \ref{AAforXY}.

\section{Wallman compactification}
\label{sectionWallman}

Throughout the whole paper, we will assume that $(T,\tau_T)$ is a \textit{Tychonoff space} i.e. it is $T_1$ and completely separates points from closed sets (there exists a function $f$ which is $0$ on the given point and $1$ on the given closed set). Moreover, we denote 

$$\Zero(T) := \bigg\{f^{-1}(0) \ : \ f \in C^b(T)\bigg\}$$

\noindent
The elements of $\Zero(T)$ are called the \textit{zero sets}. Observe that $\Zero(T)$ is closed under finite intersections. To the best of our knowledge, the importance of zero sets in the context of Wallman compactification were first noticed by Gillman and Jerison in \cite{GillmanJerison}, particularly chapter 6. The essence of $\Zero(T)$ is grasped by the following lemma, which appears (as a part of a proof) in \cite{Walker} on page 23.

\begin{lemma}
For any $A,B \in \Zero(T)$ such that $A \cap B = \emptyset$ there are $U, V \in \tau_T$ such that $A \subset U, \ B \subset V$ and $U \cap V = \emptyset$.
\label{zerosetlemma}
\end{lemma}

Intuitively, zero sets provide a substitute for the separation axiom $T_4$ (normality). On the basis of these sets, Gillman and Jerison built the Wallman compactification, which was originally carried out using all closed subsets of $T$ (comp. \cite{Wallman}).

A family $\Ufamily \subset \Zero(T)$ is called a \textit{Wallman ultrafilter} or \textit{$\omega$-ultrafilter} if

\begin{description}
	\item[\hspace{0.4cm} ($\omega$1)] Any finite intersection of elements in $\Ufamily$ is nonempty.
	\item[\hspace{0.4cm} ($\omega$2)] The family $\Ufamily$ is maximal.
\end{description}

We will denote the set of all $\omega$-ultrafilters on $T$ by $\Wall(T)$. It is easy to observe that every point $t \in T$ determines a \textit{principal $\omega$-ultrafilter} 

$$\Pfamily_t := \bigg\{ f^{-1}(0) \in \Zero(T) \ : \ f(t) = 0 \bigg\}$$

\noindent
The function $\wp : T \rightarrow \Wall(T)$ given by 

$$\forall_{t \in T} \ \wp(t) := \Pfamily_t$$

\noindent
is called a \textit{principal function}.

The first step in introducing the topology on $\Wall(T)$ is defining the operator $\ast : \tau \rightarrow 2^{\Wall(T)}$ with 

$$\forall_{U \in \tau_T} \ U^{\ast} := \bigg\{ \Ufamily \in \Wall(T) \ : \ T\backslash U \not\in \Ufamily \bigg\}$$

It is a trivial observation, yet still useful, that

\begin{gather}
\forall_{U \in \tau_T} \ t \in U \ \Longleftrightarrow \ \wp(t) \in U^{\ast}
\label{principalfilterinUstar}
\end{gather}

Among the other properties of $\ast$-operator, we recall that if $U, V \in \tau_T$ then $(U \cap V)^{\ast} = U^{\ast} \cap V^{\ast}$ and $(U \cup V)^{\ast} = U^{\ast} \cup V^{\ast}$. Furthermore, if $U \subset V$ then $U^{\ast} \subset V^{\ast}$, so the operator $\ast$ is increasing with respect to inclusion. We may conlcude that the family of sets $U^{\ast}$ where $U \in \tau_T$, is a topological base. The topology, which is thus introduced on $\Wall(T)$, is called \textit{Wallman topology} and we denote it by $\tau^{\ast}$.

\begin{thm}(Wallman compactification)\\
The pair $(\Wall(T),\wp)$ is the Stone-\v{C}ech compactification. 
\end{thm}

We omit the proof, sketching only 2 of its aspects. First, in order to establish that $(\Wall(T),\tau^{\ast})$ is Hausdorff, we use lemma \ref{zerosetlemma}. This part heavily relies on zero sets. Originally, the construction was done for normal space $T$ by Wallman in \cite{Wallman}. It was Gillman and Jerison, who realized the importance of zero sets thus requiring of $T$ only to be Tychonoff.

The proof usually goes on to show that $(\Wall(T),\tau^{\ast})$ is compact, and that $\wp$ is a homeomorphism of $T$ and a dense subspace of $\Wall(T)$. At the final stage, in order to prove that $(\Wall(T),\wp)$ is the Stone-\v{C}ech compactification, we can simply verify that every $f \in C^b(T)$ extends to $\hat{f} \in C(\Wall(T))$ in the sense $\hat{f} \circ \phi = f$. This characterization (among many others) can be found in \cite{Walker} on page 25.

For every $\Ufamily \in \Wall(T)$ we consider
	
\begin{gather}
\Ufamily_f := \bigg\{ A \subset \overline{\text{Im}(f)} \ : \ f^{-1}(A) \in \Ufamily \bigg\}
\label{defofuf}
\end{gather}
	
\noindent
It can be shown that the intersection $\bigcap \ \Ufamily_f$ is not only nonempty, but moreover consists of exactly one point. This point is called the \textit{limit of the function along $\omega$-ultrafilter} and denoted by $\lim_{\Ufamily} \ f$. Last but not least, we can define $\hat{f} : \Wall(T) \rightarrow \reals$ by
	
\begin{gather}
\forall_{\Ufamily \in \Wall(T)} \ \hat{f}(\Ufamily) := \lim_{\Ufamily} \ f \in \bigcap \ \Ufamily_f
\label{definitionoffhat}
\end{gather}

\noindent
The proof ends with the verification of the continuity of $\hat{f}$ and the identity $\hat{f}\circ \wp = f$.

\section{Arzel\`a-Ascoli theorem}
\label{sectionmainpart}

The current section will prove that $C^b(T)$ and $C(\Wall(T))$ are homeomorphic. This does not come as a surprise if one recalls the famous Gelfand-Naimark theorem, which can be found in \cite{GelfandNaimark}. The Stone-\v{C}ech compactification can manifest itself in a variety of forms:

\begin{itemize}
	\item $\beta T$, as an embedding of $T$ into a compact space $[0,1]^{C(T,[0,1])}$, or $\reals^{C^b(T,\reals)}$. 
	\item $\Wall(T)$, as a family of $\omega$-ultrafilters (precisely the construction presented above).
	\item $\Delta(C^b(T))$, as a family of nonzero algebra homomorphisms $\chi : C^b(T) \rightarrow \complex$, called \textit{characters}.
	\item $\mathcal{I}(C^b(T))$, as maximal ideals of algebra $C^b(T)$.
\end{itemize}

\noindent
The last two approaches were comprehensively studied in \cite{Simmons} in chapter 14. The advantage of $\Wall(T)$ over other forms of the Stone-\v{C}ech compactification is the simplicity of open sets. While $\beta T, \Delta(C^b(T))$ and $\mathcal{I}(C^b(T))$ all use some sort of weak$^{\ast}$ topology, Wallman compactification enjoys a very pleasant Wallman topology (described above), which is a bit easier to handle. 

Before focusing on the homeomorphism, we prove a vital property of the limit along ultrafilter. 

\begin{lemma}
If $f \in C^b(T)$ then for every $\eps > 0$ we have 

\begin{gather}
T \backslash \bigg\{t \in T \ : \ |f(t) - \lim_{\Ufamily} \ f| < \eps\bigg\} \not\in \Ufamily
\label{limufamilyf}
\end{gather}

\label{propertylimitalongultrafilter}
\end{lemma}
\begin{pro}

Suppose that (\ref{limufamilyf}) does not hold for some $\eps > 0$, i.e.

\begin{gather}
\bigg\{t \in T \ : \ |f(t) - \lim_{\Ufamily} \ f| \geq \eps\bigg\} \in \Ufamily
\label{negationoflimufamilyf}
\end{gather}

\noindent
We put $A := \complex \backslash B(\lim_{\Ufamily} \ f,\eps)$, which is obviously a closed set such that $\lim_{\Ufamily} \ f \not\in A$. Moreover, (\ref{negationoflimufamilyf}) means that 

$$f^{-1}(A) \in \Ufamily \ \stackrel{(\ref{defofuf})}{\Longrightarrow} \ A \in \Ufamily_f \ \stackrel{(\ref{definitionoffhat})}{\Longrightarrow} \ \lim_{\Ufamily} \ f \in A$$

\noindent
which is a contradiction. 
\end{pro}

\begin{thm}
The function $\Gamma : C^b(T) \rightarrow C(\Wall(T))$ given by

\begin{gather}
\forall_{f \in C^b(T)} \ \Gamma(f) := \hat{f}
\label{gelfandtransform}
\end{gather}

\noindent
is a homeomorphism.
\label{homeobc}
\end{thm}
\begin{pro}

At first, we prove the continuity of $\Gamma$. Fix $\Ufamily \in \Wall(T), \ \eps > 0$ and suppose that $d_{C^b(T)}(f,g) < \eps$. By lemma \ref{propertylimitalongultrafilter} we have
	
$$T\backslash \bigg\{t \in T \ : \ |f(t) - \lim_{\Ufamily} \ f| < \eps\bigg\} \not\in \Ufamily \ \hspace{0.5cm} \text{and} \hspace{0.5cm} \ T\backslash\bigg\{t \in T \ : \ |g(t) - \lim_{\Ufamily} \ g| < \eps\bigg\} \not\in \Ufamily$$
	
\noindent
which means
	
\begin{gather*}
\Ufamily \in \bigg\{t \in T \ : \ |f(t) - \lim_{\Ufamily} \ f| < \eps\bigg\}^{\ast} \cap \bigg\{t \in T \ : \ |g(t) - \lim_{\Ufamily} \ g| < \eps\bigg\}^{\ast} \\
= \Ufamily \in \bigg\{t \in T \ : \ |f(t) - \lim_{\Ufamily} \ f| < \eps \ \wedge \ |g(t) - \lim_{\Ufamily} \ g| < \eps\bigg\}^{\ast} 
\end{gather*}
	
\noindent
Consequently, we obtain
	
$$\Ufamily \in \bigg\{t \in T \ : \ |f(t) - g(t) - (\lim_{\Ufamily} \ f - \lim_{\Ufamily} \ g)| < 2\eps\bigg\}^{\ast}$$
	
\noindent
and conclude that 
	
$$\Ufamily \in \bigg\{t \in T \ : \ |\lim_{\Ufamily} \ f - \lim_{\Ufamily} \ g| < 2\eps + |f(t) - g(t)|\bigg\}^{\ast} \ \Longrightarrow \ \Ufamily \in \bigg\{t \in T \ : \ |\lim_{\Ufamily} \ f - \lim_{\Ufamily} \ g| < 3\eps\bigg\}^{\ast}$$

Finally, the set $\bigg\{t \in T \ : \ |\lim_{\Ufamily} \ f - \lim_{\Ufamily} \ g| < 3\eps\bigg\}$ cannot be empty, since $\Ufamily \not\in \emptyset^{\ast}$. Hence, $|\lim_{\Ufamily} \ f - \lim_{\Ufamily} \ g| < 3\eps$ and since the choice of $\Ufamily$ was arbitrary, we establish that $d_{C(\Wall(T))}(\hat{f},\hat{g}) < 3\eps$, proving continuity of $\hat{f}$.

For surjectivity, let $F \in C(\Wall(T))$. If we set $f = F \circ \wp$ then due to the fact that $\wp : T \rightarrow \wp(T)$ is a homeomorphism, we have $f \in C^b(T)$. Both functions $\Gamma(f)$ and $F$ agree on a dense set $\wp(T)$, hence by continuity they are equal $\Wall(T)$.

In order to prove that $\Gamma$ is an injection, suppose that $\Gamma(f) = \Gamma(g)$. This means that $\lim_{\Ufamily} \ f = \lim_{\Ufamily} \ g$ for every $\Ufamily \in \Wall(T)$. Focusing on the principal $\omega$-ultrafilters, we immediately conclude that $f \equiv g$. 

It remains to prove the continuity of $\Gamma^{-1}$. If we suppose that 
	
$$d_{C(\Wall(T))}(\Gamma(f),\Gamma(g)) \leq \eps \ \Longleftrightarrow \ \forall_{\Ufamily \in \Wall(T)} \ |\lim_{\Ufamily} \ f - \lim_{\Ufamily} \ g| \leq \eps$$
	
\noindent
It suffices to put $\Ufamily = \Pfamily_t$ where $t \in T$ in order to obtain $d_{C^b(T)}(f,g) \leq \eps$, which ends the proof. 
\end{pro}

The function $\Gamma$ defined in (\ref{gelfandtransform}) is actually a well-know Gelfand transform. As mentioned before, the celebrated Gelfand-Naimark theorem states that $\Gamma$ is in fact a $\ast$-isometry between $C^b(T)$ and $C(\Wall(T))$. However, we need not resort to such heavy machinery. The knowledge that $\Gamma$ is 'merely' a homeomorphism and the corollary below is sufficient for all our considerations.

\begin{cor}
If the set $\Ffamily \subset C^b(T)$ is relatively compact then $\hat{\Ffamily} := \Gamma(\Ffamily)$ is relatively compact and vice versa. 
\label{coronrlecompactness}
\end{cor}

We are ready to state and prove the culminating theorem of this paper.

\begin{thm}(Arzel\`a-Ascoli via Wallman compactification)\\
A family $\Ffamily \subset C^b(T)$ is relatively compact if and only if 

\begin{description}
	\item[\hspace{0.4cm} (AA1)] $\Ffamily$ is pointwise bounded, i.e. the set $\{f(t) \ : \ f \in \Ffamily\}$ is bounded for every $t \in T$
	\item[\hspace{0.4cm} (AA2)] $\Ffamily$ is $\omega$-equicontinuous, i.e. 
	
	$$\forall_{\substack{\Ufamily \in \Wall(T) \\ \eps > 0}} \ \exists_{\substack{V \in \tau_T \\ \Ufamily \in V^{\ast}}} \ \forall_{\substack{f \in \Ffamily \\ t \in V}} \ |f(t) - \lim_{\Ufamily} \ f| < \eps$$
\end{description}
\label{AA}
\end{thm}
\begin{pro}

Suppose that $\Ffamily \subset C^b(T)$ is relatively compact. By corollary \ref{coronrlecompactness}, this is equivalent to the relative compactness of $\hat{\Ffamily}$. By classical Arzel\`a-Ascoli theorem, we have that 

\begin{description}
	\item[\hspace{0.4cm} (AA$\omega$1) ] $\hat{\Ffamily}$ is pointwise bounded, i.e. the set $\{\hat{f}(\Ufamily) \ : \ f \in \Ffamily\}$ is bounded for every $\Ufamily \in \Wall(T)$
	\item[\hspace{0.4cm} (AA$\omega$2) ] $\hat{\Ffamily}$ is equicontinuous, i.e. 
	
		\begin{gather*}
		\forall_{\substack{\Ufamily \in \Wall(T) \\ \eps > 0}} \ \exists_{\substack{V \in \tau_T \\ \Ufamily \in V^{\ast}}} \ \forall_{\substack{f \in \Ffamily \\ \Vfamily \in V^{\ast}}} \ |\hat{f}(\Vfamily) - \hat{f}(\Ufamily)| < \eps 
		\end{gather*}
\end{description}

At first we observe that pointwise boundedness at every point, in view of equicontinuity, is equivalent to pointwise boundedness on a dense set $\wp(T)$. Thus the condition \textbf{(AA$\omega$1)} can be replaced with \textbf{(AA1)}.

Observe that the implication

\begin{gather}
\forall_{\Vfamily \in V^{\ast}} \ |\lim_{\Vfamily} \ f - \lim_{\Ufamily} \ f|< \eps \ \Longrightarrow \ \forall_{t \in V} \ |f(t) - \lim_{\Ufamily} \ f|< \eps
\label{implicationeps}
\end{gather}

\noindent
holds. Indeed, if $t \in V$ then by (\ref{principalfilterinUstar}) we know that $\Pfamily_t \in V^{\ast}$. Since $\lim_{\Pfamily_t} \ f = f(t)$, we conclude that (\ref{implicationeps}) holds and thus \textbf{(AA$\omega$2)} implies \textbf{(AA2)}.

We aim to show that 

\begin{gather}
\forall_{t \in V} \ |f(t) - \lim_{\Ufamily} \ f| < \eps \ \Longrightarrow \ \forall_{\Vfamily \in V^{\ast}} \ |\lim_{\Vfamily} \ f - \lim_{\Ufamily} \ f| \leq \eps
\label{whatweprove}
\end{gather}

\noindent
which will prove that \textbf{(AA2)} implies \textbf{(AA$\omega$2)}. We take $\Vfamily \in V^{\ast}$ and assume that $V \subset \{t \in T \ : \ |f(t) - \lim_{\Ufamily} \ f| < \eps\}$, which implies 
	
\begin{gather}
\Vfamily \in \bigg\{t \in T \ : \ |f(t) - \lim_{\Ufamily} \ f|< \eps\bigg\}^{\ast} 
\label{intersect1}
\end{gather}
		
Moreover, by lemma \ref{propertylimitalongultrafilter} we have
	
\begin{gather}
T\backslash \bigg\{t \in T \ : \ |f(t) - \lim_{\Vfamily} \ f| < \eta\bigg\} \not\in \Vfamily \ \Longleftrightarrow \ \Vfamily \in \bigg\{t \in T \ : \ |f(t) - \lim_{\Vfamily} \ f| < \eta\bigg\}^{\ast}
\label{intersect2}
\end{gather}
	
\noindent
for every $\eta > 0$. Intersecting (\ref{intersect1}) and (\ref{intersect2}) we obtain
	
$$\Vfamily \in \bigg\{t \in T \ : \ |\lim_{\Ufamily} \ f - \lim_{\Vfamily} \ f| < \eps + \eta\bigg\}^{\ast}$$
	
\noindent
for every $\eta > 0$. Reasoning as before, the set $\{t \in T \ : \ |\lim_{\Ufamily} \ f - \lim_{\Vfamily} \ f| < \eps + \eta\}$ cannot be empty, since $\Vfamily \not\in \emptyset^{\ast}$. Since the choice of $\eta$ was arbitrary, we conclude that $|\lim_{\Ufamily} \ f - \lim_{\Vfamily} \ f| \leq \eps$. We proved (\ref{whatweprove}) and consequently, \textbf{(AA$\omega$2)} is equivalent to \textbf{(AA2)}.
\end{pro}

As a final note, let us compare the obtained result with the previous work of the author together with Bogdan Przeradzki. In \cite{KrukowskiPrzeradzki}, the following generalization of Arzel\`a-Ascoli theorem was proved:

\begin{thm}(Arzel\`a-Ascoli for $\sigma$-locally compact Hausdorff space)\\
Let $(X,\tau_X)$ be $\sigma$-locally compact Hausdorff space and $(Y,d_Y)$ be a metric space. The set $\Ffamily \subset C^b(X,Y)$ is relatively compact iff

\begin{description}
	\item[\hspace{0.4cm} (KP1)] $\Ffamily$ is pointwise relatively compact
	\item[\hspace{0.4cm} (KP2)] $\Ffamily$ is equicontinuous
	\item[\hspace{0.4cm} (KP3)] $\Ffamily$ satisfies $C^b(X,Y)$-extension property, i.e. 
	
	$$\forall_{\eps > 0} \ \exists_{\substack{D \Subset X \\ \delta > 0}} \ \forall_{f,g \in \Ffamily} \ d_{C^b(D,Y)}(f,g) < \delta \ \Longrightarrow \ d_{C^b(X,Y)}(f,g) < \eps$$
	
	\noindent
	where $D \Subset X$ means that $D$ is a compact subset.
\end{description}
\label{AAforXY}
\end{thm}

In our considerations, we relaxed the assumptions on $X$ to being simply a Tychonoff space and took $Y = \complex$. We observe that \textbf{(AA2)} corresponds to \textbf{(KP2)} when we consider only principal $\omega$-ultrafilters. Indeed, we have 

\begin{gather*}
\forall_{\substack{t_{\ast} \in T \\ \eps > 0}} \ \exists_{\substack{V \in \tau_T \\ \Pfamily_{t_{\ast}} \in V^{\ast}}} \ \forall_{\substack{f \in \Ffamily \\ t \in V}} \ |f(t) - \lim_{\Pfamily_{t_{\ast}}} \ f| < \eps \\
\Longleftrightarrow \ \forall_{\substack{t_{\ast} \in T \\ \eps > 0}} \ \exists_{V_{t_{\ast}} \in \tau_T} \ \forall_{\substack{f \in \Ffamily \\ t \in V_{t_{\ast}}}} \ |f(t) - f(t_{\ast})| < \eps
\end{gather*} 

\noindent
which is exactly \textbf{(KP2)}. Consequently, the rest of $\omega$-ultrafilters in $\Wall(T)$, sometimes referred to as \textit{free $\omega$-ultrafilters}, play the same role as \textbf{(KP3)}.

\Addresses
\end{document}